\documentclass[11pt]{amsart}

\usepackage{amsmath}
\usepackage{amssymb}
\usepackage{graphicx}
\usepackage[initials]{amsrefs}
\usepackage{fancyhdr}
\usepackage{subcaption}
\usepackage{algorithm}
\usepackage{algpseudocode}
\usepackage{tikz}
\usepackage{multirow}
\usetikzlibrary{shapes.geometric, arrows.meta, positioning}

\newtheorem{theorem}{Theorem}[section]

\newtheorem{lemma}[theorem]{Lemma}

\theoremstyle{definition}
\newtheorem{definition}[theorem]{Definition}

\numberwithin{equation}{section}

\title[RO Framework for ground coverage in ASN]{Robust Optimization Framework for Ground Coverage in Aerial Sensor Networks}

\author[V. Datta]{Vanshika Datta}
\address[V. Datta]{Department of Mathematics, Indian Institute of Technology Kharagpur, Kharagpur, West Bengal, India}
\email{{\tt vanshikadutta28@gmail.com}}

\author[C. Nahak]{C. Nahak}
\address[C. Nahak]{Department of Mathematics, Indian Institute of Technology Kharagpur, Kharagpur, West Bengal, India}
\email{\tt cnahak@maths.iitkgp.ac.in}

\author[J. C. Yao]{J. C. Yao}
\address[J. C. Yao]{Research Center for Interneural Computing, China Medical University Hospital, China Medical University, Taichung 40402, Taiwan}
\email{\tt yaojc@mail.cmu.edu.tw}

\keywords{Robust optimization; Radius of robust feasibility; Voronoi diagram; Directional sensor network; Data uncertainty}

\begin{document}

\begin{abstract}
Sensors play a critical role in environmental monitoring, but their coverage performance degrades significantly under spatial uncertainty. This article proposes a robust optimization framework for maximizing ground coverage in aerial directional sensor networks subject to sensor displacement. Each aerial sensor projects a truncated sector on the ground, parameterized by its altitude, field of view, and orientation. To explicitly capture robustness against positional uncertainty, we adopt the radius of robust feasibility (RRF) as a quantitative measure of tolerance to worst-case perturbations. The RRF formulation for aerial sensor networks is embedded directly into the coverage maximization problem, ensuring feasibility under bounded uncertainty. The resulting worst-case coverage problem is nonconvex and NP-hard; therefore, a distributed greedy orientation algorithm based on Voronoi partitioning is applied to adjust sensor orientations using only local information, while directing coverage toward high-impact regions. Simulation results demonstrate that the proposed method consistently preserves robust coverage across complex terrain, varying parameters and uncertain operating conditions, highlighting its practical significance for aerial sensing applications.
\end{abstract}

\maketitle


\section{Introduction}
Sensors are devices that can perceive, record, and respond to their surrounding environment, forming the foundation of modern monitoring systems. They are widely used in applications such as airport surveillance, border monitoring, disaster response, infrastructure inspection, and traffic management \cite{li2021coverage}. As sensors have evolved from passive data collectors to active decision-making agents, wireless sensor networks (WSNs) emerged as a critical infrastructure for real-time monitoring in resource-constrained environments \cite{sangwan2015survey,ye2008robust}. Directional sensor networks (DSNs) represent a specialized class of WSNs. They use sensors such as cameras, infrared, ultrasonic, and LiDAR sensors that have limited fields of view \cite{adachi2022cooperative,li2022uav}. These sensors provide focused and high-precision coverage but require careful placement and orientation. As a result, coverage modeling in DSNs becomes more complex than in omnidirectional sensing \cite{sung2014voronoi}. The challenge is further amplified in aerial sensor networks (ASNs), where sensors mounted on UAVs operate in three-dimensional space. Their ground coverage depends on altitude, sensing range, and yaw angle, leading to fan-shaped or truncated sensing regions \cite{liang2025enhanced}.

As sensors become integral to IoT and cyber-physical systems, the quality and reliability of sensor data directly affect system performance \cite{Sobin}. IoT-enabled WSNs have transformed domains such as agriculture \cite{Dasig}, energy systems \cite{Pawar}, and healthcare \cite{Krishnamoorthy,Saeed2022}. However, sensors are often deployed in harsh or remote environments. Aging, calibration errors, improper installation, and environmental stress can introduce faults \cite{Khan}. These faults may be gradual or abrupt and can severely degrade decision-making if not detected in time \cite{Atmaji}. Consequently, robust sensor placement, fault detection observers, and learning-based fault diagnosis have become central topics in sensor network research \cite{2,3,1}. Beyond faults, uncertainty is an inherent challenge in real-world sensor deployments. Sensor location errors, environmental disturbances, and communication limitations can significantly degrade sensing performance \cite{[12]}. In IoT-based services, such uncertainties can lead to coverage gaps, misaligned sensing directions, and energy inefficiency. These issues are particularly severe in GPS-denied or cluttered environments such as forests and disaster zones. Addressing these challenges requires explicit modeling of uncertainty rather than relying on idealized assumptions \cite{wang2017coverage,[13]}. Although coverage optimization in WSNs has been widely studied \cite{sangwan2015survey}, most existing approaches assume deterministic sensor locations and orientations \cite{sung2014voronoi}. Grid-based and geometric deployment strategies aim to improve coverage efficiency \cite{mahfouz2023novel}, but they remain sensitive to small perturbations in sensor placement. For instance, adaptive orientation strategies for camera networks were proposed for airport surveillance in \cite{li2021coverage}, yet location uncertainty was not considered. This limitation motivates the need for uncertainty-aware design frameworks. An uncertainty-aware approach is therefore essential for efficient data collection and long-term sustainability in sensor networks. Designs that focus only on nominal optimality often fail when deployment conditions change.

Robust optimization provides a principled framework for addressing uncertainty by embedding it directly into the optimization model, ensuring feasibility for all realizations within predefined uncertainty sets \cite{soyster1973convex}. The framework was later extended into tractable formulations by Ben-Tal and Nemirovski \cite{ben2000robust}, and further developed for practical applications \cite{ben2009robust}. Robust optimization applies to continuous, discrete, and mixed-integer problems \cite{ben2002robust,kouvelis2013robust}, and is particularly suitable for data-sparse or high-risk environments \cite{gabrel2014recent}. Unlike stochastic programming, it does not rely on accurate probability distributions, which are often unavailable in field deployments. In sensor networks, robust optimization has primarily been applied to energy-efficient design and network longevity \cite{luo2020maximizing}, system consistency and stability under uncertainty, and energy-aware sensing strategies \cite{eren2025simulation}. Some studies have also examined robustness against constraint perturbations \cite{goberna2018guaranteeing}. However, robustness is rarely integrated directly into directional coverage models, particularly in three-dimensional sensing environments. Our earlier work introduced the radius of robust feasibility (RRF) for two-dimensional directional sensor networks using linearized formulations \cite{vanshika}; however, scalability and aerial deployment aspects remained unaddressed.

A key limitation of robust optimization is the lack of a quantitative measure that indicates how much uncertainty a solution can tolerate. RRF addresses this limitation by quantifying the maximum size of the uncertainty set for which all constraints remain satisfied. The concept was originally introduced for linear semi-infinite programs \cite{goberna2014robust} and later extended to more general convex and nonconvex settings \cite{goberna2022radius}. From a sensing perspective, RRF provides a clear robustness margin, allowing designers to certify feasibility under worst-case perturbations rather than relying on nominal configurations. Beyond its theoretical significance, RRF serves as a practical design tool for reliable systems under uncertainty by characterizing the largest admissible uncertainty region, often modeled as a norm ball or a compact convex set, within which feasibility is preserved. Several variants of RRF have been studied, including radii for multi-objective efficiency and constraint-level robustness guarantees \cite{goberna2015robust}. Such formulations are especially relevant in applications where maintaining operability under uncertainty is critical.

In distributed sensor networks, deployment uncertainty naturally arises due to terrain irregularities, environmental effects, and sensor mobility. RRF enables these uncertainties to be modeled explicitly as tolerance regions around nominal sensor locations \cite{ye2008robust}, leading to sensing strategies that remain feasible even under small spatial displacements. Recent work demonstrated that shifting from uncertain coefficients to uncertain locations represents a significant conceptual advancement \cite{ridolfi2023radius}. Methodologically, RRF has evolved substantially since its initial formulation \cite{goberna2021calculating}. Early studies focused on ball uncertainty sets and linear semi-infinite programming \cite{chuong2016robust,goberna2016radius}, while later works extended the concept to compact convex and spectrahedral uncertainty sets \cite{chen2020radius}. More recent studies introduced spatial interpretations based on uncertain Voronoi structures \cite{sung2014voronoi}. However, none of these works considered directional sensing models or their extension to three-dimensional aerial sensor networks, which motivates the present study.


While recent works have investigated RRF computation in mixed-integer or abstract optimization settings \cite{chuong2017exact,goberna2022radius,liers2022radius}, which primarily focused on linear programs or coefficient uncertainty, none have addressed aerial directional sensing models. To the best of our knowledge, this is the first work to integrate RRF, constraint-based robustness enforcement, and Voronoi-partitioned greedy orientation in aerial directional sensor networks.


In this article, we investigate the problem of coverage maximization in aerial directional sensor networks by incorporating robustness against uncertainties in sensor deployment and terrain conditions. Starting from a nominal coverage formulation based on Voronoi-partition-driven orientation strategies, we develop a robust optimization framework by integrating the RRF concept to account for bounded uncertainties in sensor locations and environmental variations. This formulation enables robustness guarantees to be embedded into the orientation-aware coverage design, treating sensor orientations as decision variables. The resulting coverage optimization problem is nonconvex and computationally challenging. To address this, we explore a distributed greedy orientation algorithm that utilizes only local neighborhood information, thereby ensuring scalability and practical applicability in large-scale aerial sensor deployments. The effectiveness of the proposed framework is demonstrated through extensive simulation studies under both nominal and perturbed deployment scenarios. The proposed algorithm outperforms the IDS-based and random orientation methods, delivering higher and more reliable coverage even under worst-case sensor displacements. Although recent studies have explored RRF, mainly in the context of linear programs, mixed-integer formulations, or abstract coefficient-uncertainty settings, their application to aerial directional sensing models has not been previously addressed. To the best of our knowledge, this work presents the first unified framework that combines RRF-based robustness enforcement and a Voronoi-partitioned greedy orientation algorithm for coverage in aerial sensor networks.


The remainder of the article is organized as follows. Section 2 introduces the notation and preliminary results. The modeling framework and problem setup are described in Section 3. Section 4 introduces the orientation adjustment strategy and the distributed optimization algorithm. Section 5 provides simulation results and performance evaluation. Concluding remarks of the article and future work are presented in Section 6.

\section{Preliminaries}  \label{sec:2} 
This section introduces the formal notation used throughout the article, along with the Voronoi-based spatial decomposition, sensing geometry, and the adopted sensing model.

\subsection{Notations and basic results}
Let $\mathbb{R}^n$ denote the $n$-dimensional Euclidean space, $\mathbb{B}_n$ the open unit ball in $\mathbb{R}^n$, and $||\cdot||$ the Euclidean norm. For any vector $a$, its transpose is denoted by $a^T$. The Euclidean distance between two points $x,y \in \mathbb{R}^n$ is given by $d(x,y)=||x-y||$. For a set $S \subset \mathbb{R}^n$, we denote its closure, convex hull, and conical hull by $cl(S)$, $conv(S)$, and $cone(S)$, respectively. The set of nonnegative functions defined on $S$ is denoted by $\mathbb{R}_+^S$, while $\mathbb{R}_+^{(S)}$ represents the set of finitely supported nonnegative functions on $S$.\\
A parametric linear system in the face of data uncertainty in its constraints, denoted by $\sigma^{\alpha }  $, can be captured as follows: 
\begin{equation}
    \sigma ^{\alpha } := \left\{ a^T_ix \leq b_i \right\}; \;(a_i,b_i) \in \mathcal{U} _i ^{\alpha} \subset \mathbb{R} ^{n+1}; \;i=1:p \;,
\end{equation}
where $(a_i,b_i) $ for $i=1:p$ are uncertain vectors and $\mathcal{U} _i ^{\alpha};\; i=1:p$ are the uncertain regions, which are typically assumed to be compact and convex.\\
The robust counterpart (deterministic form) of the corresponding system $\sigma ^{\alpha} $ with $U_i^{\alpha}=(\bar{ a_i}, \bar{b_i} ) + \alpha B_{n+1}, \; i=1:p$ and the corresponding feasible set $F_R^ {\alpha}$ for some $\alpha \in \mathbb{R} $ is:
\begin{equation}
    \sigma _ R ^ {\alpha} := \left\{ a_i^Tx \leq b_i ,\;  \forall (a_i,b_i) \in U_i^{\alpha},\; i=1:p \right\}.
\end{equation}
The radius of robust feasibility $\rho$, of an uncertain system $\sigma ^{\alpha}$, determines the maximum level of perturbation that a system can tolerate while remaining feasible. It is formally defined as follows:
\begin{definition}(Radius of robust feasibility)\label{rrf} Consider a parametric linear system in face of data uncertainty $\sigma ^{\alpha} $. Let $F_R^{\alpha}$ denote the feasible set of the robust counterpart of $\sigma ^{\alpha}$. Then, the RRF for the system is given by: 
\begin{equation}
     \rho = \sup \{ \alpha \in \mathbb{R}_+ : (F_R ^ \alpha) \text{ is nonempty} \}.
\end{equation} 
\end{definition}
This quantifies the maximum perturbation level under which the system remains robustly feasible. A fundamental concept in computing RRF is the Minkowski function, defined as follows:
\begin{definition}(Minkowski function)
    Let $\omega \subset \mathbb{R} ^ n $ be a convex set containing $0_n$ in its interior. Then the Minkowski or gauge function of $\omega$ denoted by $\phi_{\omega}$, where  $\phi_{\omega} : \mathbb{R} ^n \to \mathbb{R} _+ := [0, + \infty [ $ is given by: 
    \begin{equation}
        \phi_{\omega} (x) := \inf{ \left\{ t>0 : x \in t \omega \right\} },\; x \in \mathbb{R} ^n.
    \end{equation}
\end{definition}
The following lemma provides some properties of the Minkowski function.
\begin{lemma}[See \cite{chuong2017exact}, Lemma 2.2] Let $\omega \subset \mathbb{R} ^n $ be a convex set such that its interior contains $0_n$, then the following properties hold:
\begin{itemize}
    \item[$(1)$] $\phi_{ \omega}$ is sublinear and continuous.
    \item[$(2)$] $\left\{ x \in \mathbb{R} ^n : \phi _{\omega} (x) \leq 1 \right\} = cl (\omega)$, where $cl (\omega ) $ stands for the closure of $\omega$.
    \item[$(3)$] If in addition, $\omega$ is bounded and symmetric, then $\phi _ {\omega }: = || \cdot ||$ is a norm on $\mathbb{R} ^n$ generated by $\omega$.
\end{itemize} 
\end{lemma}
Exact formula of RRF for an arbitrary set $S \subset \mathbb{R}^{n+1}$ has been obtained in the literature in terms of the so-called hypographical set $H(\bar a, \bar b)$ of the nominal system $\sigma$: 
\begin{equation}
    H(\bar a, \bar b) := conv \left\{ (- \bar{a_i} , -\bar{b_i}): i=1:p \right\} + \mathbb{R} _+ (0_n, -1),
\end{equation}
where $\bar a = ( \bar{a_1},\bar{a_2}, \cdots, \bar{a_p} ) \in (\mathbb{R} ^n)^p \text{ and } \bar b = ( \bar{b_1},\bar{b_2}, \cdots, \bar{b_p} ) \in \mathbb{R} ^p$.
\begin{theorem} 
    If the nominal system is feasible, the RRF of $\sigma^{\alpha}$ is:
    \begin{equation}
        \rho = \inf_{(a,b) \in H(\bar{a},\bar{b})} \phi_Z (a,b),
    \end{equation}
    where $Z \subset \mathbb{R} ^n $ is a convex and compact set containing $0_n$ in its interior, see \cite{chuong2017exact}. 
\end{theorem}
\begin{figure}
    \centering
    \includegraphics[width=0.4\textwidth]{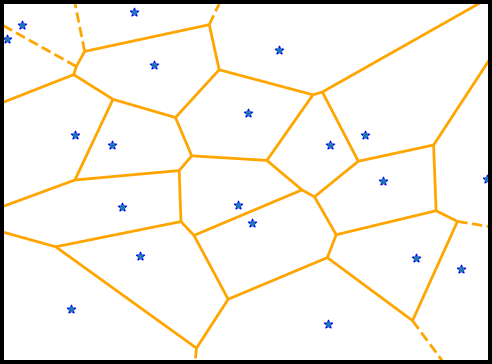}
    \caption{Voronoi diagram for a set of generators}
    \label{fig:voronoidiag}
\end{figure}

\begin{definition}(Voronoi cell) \label{voronoi_cell_def}
    Given a set of generator points $S=\{s_1, s_2, \dots, s_n\}\in \mathbb{R}^2$, the Voronoi cell $V_S(s_i)$ associated with sensor $s_i \in S$ is defined as:
    \begin{equation}
         V_S(s_i) = \{ p \in \mathbb{R}^2 : d(p, s_i) \leq d(p, s_j) \text{ for all } j \neq i \}, 
    \end{equation}
    where  $j=1:n$ and $\;d(a,b)$ is the Euclidean distance between $a$ and $b$.
\end{definition}
A Voronoi diagram is a spatial decomposition of a region based on proximity to a given set of points, called generators (sensor nodes in our case). It consists of the union of all Voronoi cells (see Figure~\ref{fig:voronoidiag}). 

\subsection{Sensor geometry and sensing model}\label{3.2}
\begin{figure}
\begin{center}
    \begin{subfigure}{.5\textwidth}
  \centering
  \includegraphics[width=.8\linewidth]{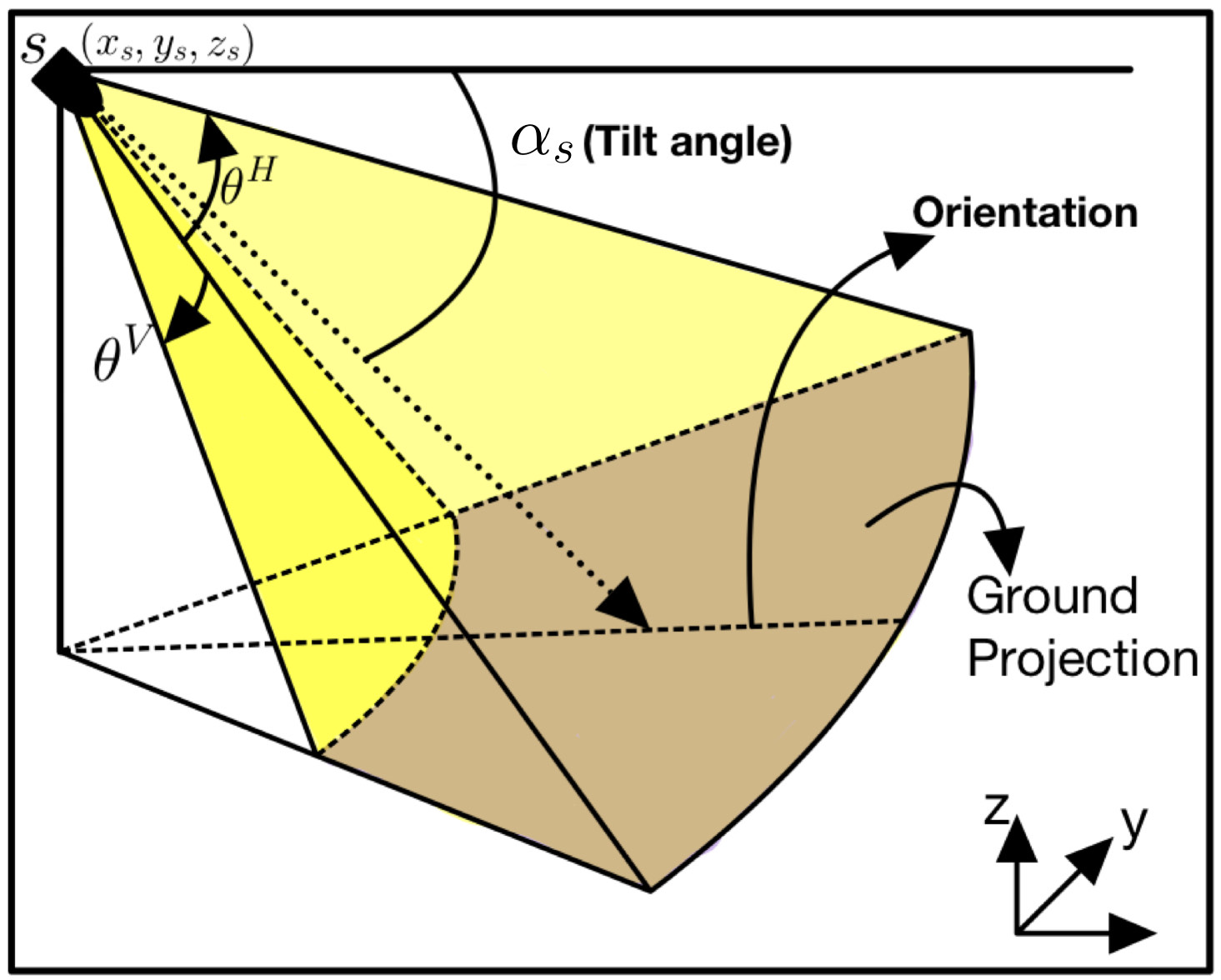}
  \caption{A 3D sensor model}
  \label{fig:3D}
\end{subfigure}%
\begin{subfigure}{.5\textwidth}
  \centering
  \includegraphics[width=.8\linewidth]{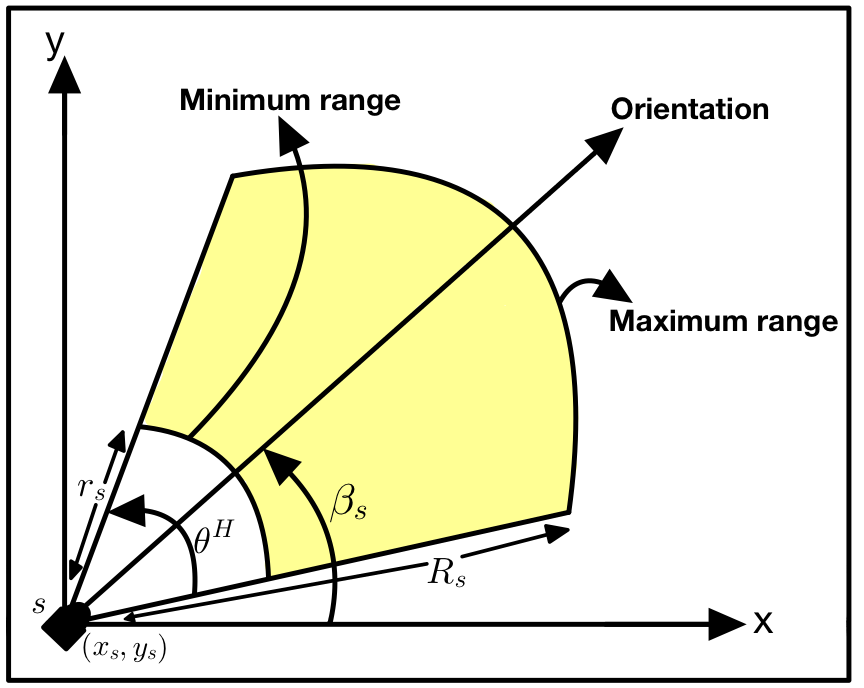}
  \caption{Ground projection of a 3D sensor}
  \label{fig:2D}
\end{subfigure}
\caption{Sensing model in 3D and its ground projection}
\label{fig: sensing model}
\end{center}
\end{figure}

Sensors are deployed in three-dimensional space at locations $(x_s, y_s, z_s) \in \mathbb{R}^3$, typically elevated to expand the field of view and mitigate occlusion effects, resulting in a top-down orientation as illustrated in Figure~\ref{fig:3D}. Each sensor is characterized by a horizontal and a vertical angle of view, denoted by $\theta_H$ and $\theta_V$, respectively. The tilt angle $\alpha_s$ is defined as the angle between the sensor’s optical axis and the $xy$-plane and controls vertical inclination, while the yaw angle $\beta_s$ specifies the horizontal orientation with respect to the positive $x$-axis, with $\beta_s = 0$ indicating alignment along the $x$ direction.

For a sensor located at $(x_s, y_s, z_s)$ with tilt angle $\alpha_s$, the three-dimensional sensing region forms a truncated pyramidal frustum, as illustrated in Figure~\ref{fig: sensing model}. Since coverage is evaluated on the ground, this volume is orthogonally projected onto the $xy$-plane, producing a truncated sector or fan-shaped footprint, shown in Figure~\ref{fig:2D}. After projection, each aerial directional sensor is represented in two dimensions by its planar location $(x_s, y_s)$, orientation $\beta_s$, horizontal field of view $\theta_H$, and sensing limits. The inner and outer sensing radii are denoted by $r_s$ and $R_s$, respectively. The resulting ground coverage region of the sensor $s$ is an annular sector $C_s(\beta_s)$, defined as:
\begin{equation}\label{sensing_area}
\begin{split}
    C_s(\beta_s)=&\left\{(x,y)\in\mathbb{R}^2 \,\middle|\,
     r_s \le \|(x,y)-p_s\| \le R_s,\right.\\
     & \;\;\;\;\;\;\;\;\;\left.\operatorname{atan2}(y-y_s,\,x-x_s)\in
     \left[\beta_s-\frac{\theta_H}{2},\,\beta_s+\frac{\theta_H}{2}\right]\right\},
\end{split}
\end{equation}
where $\operatorname{atan2}(y - y_s, x - x_s)$ denotes the angle between the vector from the sensor $s$ to the point $(x,y)$ and the positive $x$-axis.

Let $T \subset \mathbb{R}^3$ denote the location set of aerial sensors, with $|T| \geq 2$, and let $S$ be its projection onto the ground plane. We henceforth consider only the sensor’s planar coordinates. Each sensor location in $S$ is denoted by $s = (x_s, y_s) \in \mathbb{R}^2$. A point $p = (x, y) \in \mathbb{R}^2$ on the ground is said to be covered by sensor $s$ if it satisfies both the radial and angular constraints of $A_s(\beta_s)$. Hence, the effective sensing region of each sensor is fully determined by the parameters $(r_s, R_s, \theta_H, \beta_s)$. Specifically, sensor $s$ covers point $p$ if the following two conditions are met:
\begin{itemize}
    \item[$(1)$] The Euclidean distance from $p$ to the sensor lies within the sensor's detection range:
    \begin{equation}
        r_s \leq \sqrt{(x - x_s)^2 + (y - y_s)^2} \leq R_s.
    \end{equation}
    \item[$(2)$] The angle between the sensor's forward orientation $\beta_s$ and the direction vector $\vec{sp}$ is within the sensor’s angular field of view:
    \begin{equation}
        \left| \operatorname{atan2}(y - y_s, x - x_s) - \beta_s \right| \leq \frac{\theta_H}{2}.
    \end{equation}
\end{itemize}

\subsection{Problem Overview and Assumptions}
We study the problem of maximizing ground coverage in an aerial directional sensor network deployed over a region of interest (ROI), where each sensor observes a sector-shaped sensing region and its ground-projected location may deviate from the nominal position due to environmental disturbances or deployment inaccuracies. The objective is to determine the sensing orientations of the deployed sensors so that the total effectively covered ground area is maximized while maintaining feasibility under bounded location uncertainty.\\
To organize sensing responsibility and reduce redundant overlap, the ROI is partitioned using the Voronoi diagram constructed from the nominal sensor locations, assigning each sensor a spatial subregion for coverage evaluation. Although Voronoi partitioning does not necessarily produce the globally optimal assignment of sensing regions, it provides a simple and computationally tractable framework to associate each sensor with a local spatial responsibility region and to limit redundant overlap in coverage analysis. The convexity of Voronoi cells implies that the farthest point from a sensor within its cell lies at one of the cell vertices (Theorem 3.2 in \cite{vanshika}), which allows the orientation analysis to be performed using only these vertices. This vertex-based characterization significantly simplifies the coverage evaluation while preserving tractability. Moreover, since the Voronoi partition assigns every sensor a nonempty spatial cell, the nominal configuration does not produce infeasible regions, ensuring that coverage analysis and optimization can be carried out independently for each sensor within its associated cell.

\section{Proposed Methodology}\label{sec:3}

This section presents the theoretical modeling and optimization framework for robust coverage in aerial directional sensor networks under uncertainty and outlines the problem formulation and solution flow. Table~\ref{tab:notation} explains the principal notation used throughout the paper.

\begin{table}
\centering
\caption{Table of notation.}
\begin{tabular}{|p{3cm}|p{9cm}|}
\hline
\textbf{Symbol} & \textbf{Description} \\
\hline

$\mathbb{R}^n$ & $n$-dimensional Euclidean space. \\

\hline

$\mathbb{B}_n$ & Unit ball in $\mathbb{R}^n$ with respect to the Euclidean norm. \\

\hline

$||\cdot||$ & Euclidean norm. \\

\hline

$d(x,y)$ & Euclidean distance between points $x$ and $y$. \\

\hline

ROI & Region of interest where coverage is evaluated. \\

\hline

$m$ & Total number of sensors in the deployment. \\

\hline

$\phi_{\omega}$ & Minkowski function associated with a convex set $\omega$. \\

\hline

$s$ & A particular sensor for reference. \\

\hline

$S = \left\{ s^o_1, \cdots, s^o_m \right\}$ & Nominal location set of ground-projected sensors. \\

\hline

$s^o_i=(x^o_i,y^o_i)$ & Nominal location of the $i^{th}$ sensor in $\mathbb{R} ^2$. \\

\hline

$s^w_i$ & Worst case location of  the $i^{th}$ sensor in $\mathbb{R} ^2$. \\

\hline

$s^{\rho}_i$ & Robust case location of the $i^{th}$ sensor in $\mathbb{R} ^2$. \\

\hline

$o_i$ & Orientation of the $i^{th}$ sensor after ground projection. \\

\hline

$V_i$ & Selected vertex of orientation of the $i^{th}$ sensor. \\

\hline

$A_i(s_i,o_i)$ & Set of all covered points by the $i^{th}$ sensor in $\mathbb{R} ^2$. \\

\hline

$r_s$ & Minimum sensing range of a particular sensor $s$. \\

\hline

$R_s$ & Maximum sensing range of a particular sensor $s$. \\

\hline

$VC(s^o_i)$ & Voronoi cell of the $i^{th}$ sensor corresponding to the nominally located sensors. \\

\hline

$\mathcal{U}_i^{\alpha}$ & Uncertainty set describing admissible perturbations of the $i^{th}$ sensor with radius $\alpha$. \\

\hline

$\alpha$ & Radius of the uncertainty set. \\

\hline

$\rho_{i}$ & Radius of Robust Feasibility (RRF) associated with the $i^{th}$ sensor. \\

\hline

$C_i{(\beta_i)}$ & Ground coverage of the $i^{th}$ sensor. \\

\hline

\end{tabular}
\label{tab:notation}
\end{table}

\subsection{Model Construction} 

We consider a set of $m$ sensors with nominal location set $S = \left\{ s^o_1,s^o_2, \cdots, s^o_m \right\}\subset \mathbb{R} ^2$, where $s^o_i=(x^o_i,y^o_i)$ denotes the ground-projected nominal location of the $i^{th}$ sensor for $i=1:m$.

For the $i^{th}$ sensor, let $V_i$ denote the selected vertex toward which the sensor is oriented and $A_i$ denote the set of all covered points by the $i^{th}$ sensor in $\mathbb{R} ^2$, which depends on the sensor location $s_i$ and its orientation $o_i$ toward $V_i$, for $i=1:m$. The corresponding area covered by $i^{th}$ sensor is denoted by $\mathrm{Area}(A_i(s_i,o_i))$. Based on these definitions, the coverage region of each sensor must satisfy the following feasibility constraints:
\begin{align} \label{constraints}
\begin{cases}
    A_i \subseteq  VC(s^o_i):= &\left\{p \in \mathbb{R} ^2 : ||p-s^o_i|| \leq ||p-t||, \; \forall \; t\in S \right\} \\
&\hspace{5cm} \text{(Voronoi cell constraint)},\\
A_i \subseteq  C_{i}(\beta_{i}) & \hspace{4.6cm} \text{(Sensing region constraint)},\\
A_i \subseteq  \text{ROI} & \hspace{5.4cm} \text{(Boundary constraint)},
\end{cases}
\end{align}
where $C_{i}(\beta_{i})$ is defined in equation~(\ref{sensing_area}) and ROI denotes the region of interest.

\subsubsection{Nominal Coverage Model} \label{nominal_model}  

We first consider the nominal coverage optimization problem in which sensor locations are assumed to be known exactly. The nominal objective is
\begin{equation} \label{nominal_model}
\text{(Nominal Model):  } \max_{o_i} \sum_{i=1}^{m} \mathrm{Area}(A_i({s^o_i, o_i})),
\end{equation}
subject to constraints \eqref{constraints}. This formulation maximizes the total coverage of the sensor set using the Voronoi partition while preventing out-of-boundary sensing. This model is very optimistic as it assumes exact sensor locations and ignores deployment uncertainty. Therefore, it does not guarantee robust feasibility.

\subsubsection{Robust Coverage under Location Uncertainty}  \label{RC} 

To capture deployment inaccuracies, each sensor location is allowed to deviate from its nominal position ${s^o_i}$ within the uncertainty set:
\begin{align}\label{uncertainity}
    \mathcal{U}_i^{\alpha} = {s^o_i} + \alpha \mathbb{B}_2,
\end{align}
where the scalar $\alpha>0$ represents the radius of the uncertainty set and quantifies the maximum admissible deviation of the actual sensor location from its nominal position. Hence, the uncertainty set considered in our model represents perturbations in the sensor’s nominal position. Such perturbations naturally capture lateral displacements in the horizontal plane (and more generally spatial position errors). Therefore, lateral displacement uncertainty is implicitly incorporated in the model. Other possible uncertainties, such as sensing range fluctuations or altitude variations, are not explicitly considered here in order to keep the focus on location perturbations and orientation optimization.\\
In this model, to ensure robustness, we evaluate sensor coverage under the worst-case realization of location uncertainty. For each sensor location $s^o_i \in S$, we compute the minimum achievable coverage area $\mathrm{Area}(A_i)$ over all admissible locations within its uncertainty set $\mathcal{U}_i^{\alpha}$. Since $\mathcal{U}_i^{\alpha}$ is convex, the minimum coverage is attained at an extreme point of the set, that is, on its boundary. Throughout the paper, Voronoi cells are constructed using the nominal sensor locations. These cells serve as a fixed spatial reference to evaluate how coverage changes under perturbations. This observation leads to the following lemma.
\begin{lemma}[Coverage Under Uncertainty]
For a convex uncertainty set $\mathcal{U}_i^{\alpha}$, the minimum coverage of a sensor $s^o_i \in S$ is achieved at the boundary of $\mathcal{U}_i^{\alpha}$.
\end{lemma}

Accordingly, we consider the worst-case approximate location $s^w_{i}$ of the $i^{th}$ sensor corresponding to its nominal position $s^o_i$ along its sensing direction $\vec{u}$ as
\begin{equation}
    s^w_{i} = {s^o_i} + \alpha \frac{\vec u}{|\vec u|},
\end{equation}
which provides a tractable representation for robust coverage analysis. 

Therefore, the robust counterpart of the nominal model \eqref{nominal_model} is formulated. Under this model, coverage must remain feasible for all realizations $s_i \in \mathcal{U}_i^{\alpha}$. Evaluating coverage at the boundary realization $s^w_{i}$ yields the following deterministic robust counterpart:
\begin{equation} \label{RC_model}
\text{(RC Model):  }\max_{o_i} \sum_{i=1}^{m} \mathrm{Area}(A_i(s^w_{i},o_i)),
\end{equation}
subject to constraints \eqref{constraints}, which guarantees feasibility under bounded location perturbations. Although this formulation ensures robustness, it is conservative because it always assumes the largest possible perturbation for every sensor and does not quantify how much uncertainty the obtained solution can actually tolerate. 

\subsubsection{RRF-Based Robust Optimization Formulation} 

To quantify this tolerance, we employ the concept of the radius of robust feasibility (RRF), originally introduced in \cite{vanshika}. In this work, the RRF of a sensor denotes the largest perturbation radius for which the coverage constraints \eqref{constraints} remain satisfied while maintaining a prescribed coverage threshold $\delta$.

Let $\rho_{s_i}$ denote the RRF corresponding to the $i^{th}$ sensor. This quantity characterizes the maximum admissible displacement from the nominal sensor location for which the required sensing coverage remains feasible. Based on the computed RRF, we consider the following robust approximate location of the $i^{th}$ sensor along its sensing direction $\vec{u}$:
\begin{equation}
    s^{\rho}_{i} = {s^o_i} + \rho_i \frac{\vec u}{|\vec u|},
\end{equation}
which provides a tractable representation for evaluating coverage under admissible perturbations. 

Building on this idea, we introduce a robust optimization formulation that incorporates the RRF into the coverage model. Let $\rho_{\min}$ denote a user-defined robustness requirement representing the minimum perturbation level that must be tolerated by the system. Enforcing this condition ensures that the obtained sensing configuration remains feasible under location deviations of at least $\rho_{\min}$. 

The resulting optimization problem is formulated as
\begin{equation} \label{rrf-rc_model}
\text{(RRF-RC Model):  }\max_{o_i} \sum_{i=1}^{m} \mathrm{Area}(A_i(s^{\rho}_{i},o_i)).
\end{equation}
subject to constraints \eqref{constraints}. 

\begin{figure}
\centering
\resizebox{0.42\textwidth}{!}{
\begin{tikzpicture}[
    node distance=1.3cm and 2.4cm,
    every node/.style={font=\scriptsize},
    arrow/.style={->,>=stealth,thin}
]

\tikzstyle{startstop} = [
    ellipse,
    minimum width=2.6cm,
    minimum height=0.8cm,
    text centered,
    draw=black
]

\tikzstyle{inputoutput} = [
    trapezium,
    trapezium left angle=70,
    trapezium right angle=110,
    minimum width=2.9cm,
    minimum height=0.8cm,
    text centered,
    draw=black
]

\tikzstyle{process} = [
    rectangle,
    minimum width=2.9cm,
    minimum height=0.8cm,
    text centered,
    draw=black
]

\tikzstyle{decision} = [
    diamond,
    aspect=2.2,
    minimum width=3.0cm,
    text centered,
    draw=black,
    inner sep=1pt
]

\node (start) [startstop] {Start};
\node (input) [inputoutput, below of=start] {Input sensor data};
\node (rrf) [process, below of=input] {Compute RRF};
\node (check) [decision, below of=rrf, yshift=-0.15cm]
    {$\rho_{min} < \rho < \rho_{max}$?};

\node (accept) [process, right of=check, xshift=3.2cm]
    {Sandwich $\rho$ within range};

\node (location) [process, below of=check, yshift=-0.3cm] {Update sensor locations};

\node (Model_solve) [process, below of=location]
    {Solve \eqref{rrf-rc_model}};

\node (output) [inputoutput, below of=Model_solve]
    {Final orientations};

\node (end) [startstop, below of=output] {End};

\draw [arrow] (start) -- (input);
\draw [arrow] (input) -- (rrf);
\draw [arrow] (rrf) -- (check);

\draw [arrow] (check.east) -- node[above, font=\scriptsize] {No} (accept.west);
\draw [arrow] (check.south) -- node[right, font=\scriptsize] {Yes} (location);

\draw [arrow] (location) -- (Model_solve);
\draw [arrow] (Model_solve) -- (output);
\draw [arrow] (accept) |- (location);
\draw [arrow] (output) -- (end);

\end{tikzpicture}
}
\caption{RRF-based robust sensor orientation framework}
\label{fig:rrf_flowchart}
\end{figure}
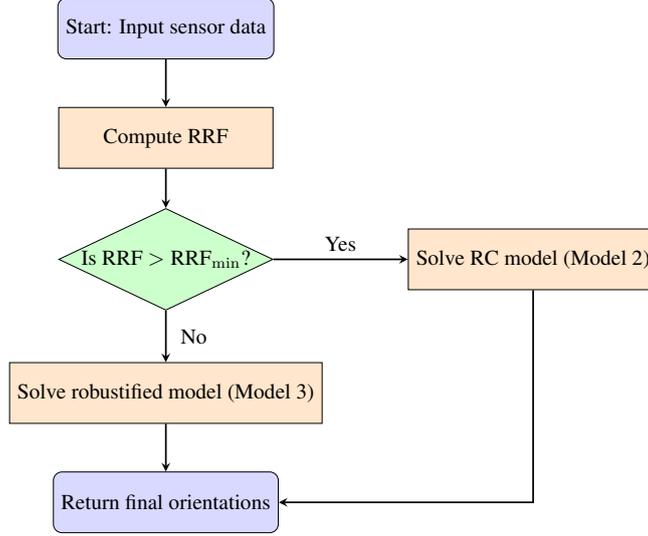

This model optimizes sensor orientations while guaranteeing a prescribed level of robustness against sensor location perturbations. The perturbation radius $\rho$ provides a unified interpretation of the three formulations presented in this section. At $\rho = 0$, the perturbed sensor location coincides with the nominal location and the formulation reduces to the nominal coverage model. In contrast, $\rho = \rho_{\max}$ corresponds to the largest admissible perturbation derived from deployment considerations, in which case the formulation reduces to the deterministic robust counterpart (RC Model), which protects against the worst-case sensor displacement but may be overly conservative. The RRF-based formulation, therefore, provides an intermediate design framework in which the effective perturbation lies within the interval
$$
\rho_{\min} \leq \rho \leq \rho_{\max},
$$
where $\rho_{\min}$ denotes the user-specified minimum robustness requirement and $\rho_{\max}$ represents the maximum admissible perturbation determined from practical considerations. This allows the sensing configuration to balance robustness and coverage performance without assuming an excessively pessimistic worst-case scenario.\\
Note that this objective maximizes the total covered area within the ROI under a uniform importance assumption. In practical scenarios, different subregions may have varying priorities or may require redundant coverage. The proposed framework can be extended to incorporate such requirements by introducing spatial weighting functions or multi-coverage constraints. The overall workflow of the proposed RRF-based robust coverage framework is summarized in Fig.~\ref{fig:rrf_flowchart}.

\section{Orientation adjustment strategy for coverage enhancement}  \label{sec:4} 
This section focuses on adjusting sensor orientations to enhance coverage while respecting the robustness guarantees established earlier. Following the Voronoi-based orientation strategy of Sung et al.~\cite{sung2014voronoi}, each directional sensor adapts its viewing direction toward critical regions of its Voronoi cell. Coverage conditions are evaluated using the sensor’s worst-case location derived from its RRF. For the $i^{\text{th}}$ sensor with location $s_i,; i=1:m,$ with nominal position $(x^o_i,y^o_i)$ and RRF $\rho_i$, oriented toward a Voronoi vertex $V_i=(a,b)$, the direction vector is $\vec{\mu}=(a-x^o_i,b-y^o_i)$, and the corresponding robust-case location $s^{\rho}_{i}=(x^{\rho}_i,y^{\rho}_i)$ is obtained by normalizing this direction and scaling it by $\rho_{i}$ as:
\begin{equation} \label{4.1}
    \begin{split}
        (x^{\rho}_i,y^{\rho}_i) &= (x^o_i,y^o_i) + \rho_{i} \frac{\vec{\mu}}{|\vec{\mu}|} \\
    &= (x^o_i,y^o_i) + \rho_{i} \left(\frac{a-x^o_i}{\sqrt{(a-x^o_i)^2 + (b-y^o_i)^2}},\frac{b-y^o_i}{\sqrt{(a-x^o_i)^2 + (b-y^o_i)^2}}\right).
    \end{split}
\end{equation}
Accordingly, coverage is evaluated using the sensor’s worst-case location and orientation. The Voronoi partition assigns primary responsibility to the nearest sensor to ensure consistent orientation decisions, while any additional coverage by neighboring sensors is naturally beneficial and does not compromise the framework. A point $p(a_0,b_0)$ is covered by the $i^{th}$ sensor if $r_s \leq d(s_{i}, p) \leq R_s$ and $\omega \leq \frac{\theta_H}{2}$, where $r_s$, $R_s$, and $\theta_H$ are the minimum range, maximum range, and angular coverage of $s_j$, and $\omega$ is the angle between $\vec{sp}$ and $\hat{\mu}$. If $o_i$ is the angle between $\hat{\mu}$ and the positive $x$-axis, $-\pi \leq o_i \leq \pi$, the second condition yields:
\begin{equation}\label{cond2}
    \begin{split}
        \vec{sp} \cdot \hat{\mu} &= ||\vec{sp} ||\; ||\hat{\mu}|| \; \cos{\omega} \\
    \implies (a_o-x^{\rho}_i) \cos{o_i}+(b_0-y^{\rho}_i) \sin{o_i}&=\sqrt{(a_0-x^{\rho}_i)^2 + (b_0-y^{\rho}_i)^2} \cos{\omega} \\   \implies (a_o-x^{\rho}_i) \cos{o_i}+(b_0-y^{\rho}_i) \sin{o_i}&\geq \sqrt{(a_0-x^{\rho}_i)^2 + (b_0-y^{\rho}_i)^2} \cos{\frac{\theta_H}{2}}.
    \end{split}
\end{equation}
Thus, under the robust case, $p(a_0,b_0)$ is covered if it satisfies $r_s \leq d(s_{i}, p) \leq R_s$ and~(\ref{cond2}).

\subsection{Classification of Voronoi-projection intersection cases}
For each sensor $s \in S$, coverage optimization is performed within its Voronoi cell by evaluating the covered area at its robust location for orientations toward all the candidate Voronoi vertices. The sensor orientation that yields the average maximum feasible coverage is selected. The sensing region is bounded by the inner and outer range arcs, denoted by $s_{in}$ and $s_{out}$, and by the sector sidelines. Depending on the intersections among Voronoi edges, sensing arcs, and sector boundaries, seven distinct geometric cases arise.
\begin{enumerate}
    \item No Voronoi-edge intersection with either arc or the sidelines.
    \item Voronoi-edge intersection with $s_{out}$ and a sideline.
    \item Voronoi-edge intersection with $s_{in}$ and a sideline.
    \item Voronoi-edge intersection with only the sidelines from inside the sense range.
    \item Voronoi-edge intersection with only the sidelines from outside the sense range.
    \item Voronoi-edge intersection with both arcs $s_{in} $ and $s_{out}$ as well as the sidelines.
    \item No intersection with any boundary element.
\end{enumerate}

\subsection{Analytical computation of covered area for each case}
Under the configuration where no Voronoi edge intersects the sensing arcs or sector sidelines, the coverage of the $i^{th}$ sensor located at $s_i$ within its Voronoi cell is computed using the sector-area expression
\begin{equation}
    A_i= \frac{\theta_H}{2} {r'_s}^2,
\end{equation}
where $\theta_H$ is the field of view and $r_s' = R_s - r_s$ denotes the effective sensing range of $i^{th}$ sensor.

For all remaining configurations, the coverage area is evaluated by decomposing the sensing region of the $i^{th}$ sensor within its Voronoi cell into simpler geometric components. Let $N$ denote the number of Voronoi vertices intersecting the sensing sector. Line segments are drawn from the robust sensor location $s^{\rho}_{i}$ to all intersection points between the sector boundaries (arcs or sidelines) and Voronoi edges, resulting in a partition into $N+1$ elementary regions, each being a triangle, a sector, or a truncated sector.

The area of each triangular component is computed using Heron’s formula
\begin{equation} \label{eqn2}
\text{Area} = \sqrt{d(d-e_1)(d-e_2)(d-e_3)},
\end{equation}
where $d=\frac{e_1+e_2+e_3}{2}$ is the semiperimeter and $e_j,; j=1:3,$ are the side lengths of the triangle.

The Voronoi edges and vertices are directly obtained from the Voronoi diagram, while the intersection points with the sensing sector are computed using the sector sideline equations. These sidelines correspond to the equality case of condition~(\ref{cond2}). Accordingly, for the $i^{th}$ sensor with robust location $s^{\rho}_{i}=(x^{\rho}_i,y^{\rho}_i)$, orientation $o_i$, angular view $\theta_H$, and sensing ranges $r_s$ and $R_s$, the sector sidelines are defined by
\begin{equation} \label{eqn4}
\begin{cases}
    (a_o-x^{\rho}_i) \cos{o_i}+(b_0-y^{\rho}_i) \sin{o_i}= \sqrt{(a_0-x^{\rho}_i)^2 + (b_0-y^{\rho}_i)^2} \cos{\frac{\theta_H}{2}}, \\
    r_s \leq d(s^{\rho}_i,p) \leq R_s.
\end{cases}
\end{equation}
Using these intersection points together with the known Voronoi geometry, the coverage area corresponding to the robust sensor position is then computed via the associated sector arc equations:
\begin{equation} \label{eqn5}
    \begin{cases}
        d(s^{\rho}_i,p) = r_s,\\
        a_0 \in [r_s \cos{\left( o_i + \frac{\theta_H}{2} \right)}+ x^{\rho}_i, r_s \cos{\left( o_i - \frac{\theta_H}{2} \right)} +x^{\rho}_i],\\
        b_0 \in [r_s \sin{\left( o_i - \frac{\theta_H}{2} \right)}+ y^{\rho}_i, r_s \sin{\left( o_i + \frac{\theta_H}{2} \right)} +y^{\rho}_i]
    \end{cases}
\end{equation}
and 
\begin{equation}
    \begin{cases}
        d(s^{\rho}_i,p) = R_s,\\
        a_0 \in [r_s \cos{\left( o_i + \frac{\theta_H}{2} \right)}+ x^{\rho}_i, r_s \cos{\left( o_i - \frac{\theta_H}{2} \right)} +x^{\rho}_i],\\
        b_0 \in [r_s \sin{\left( o_i - \frac{\theta_H}{2} \right)}+ y^{\rho}_i, r_s \sin{\left( o_i + \frac{\theta_H}{2} \right)} +y^{\rho}_i],
    \end{cases}
\end{equation}
where $p=(a_0,b_0)$ denotes a point on the arc boundary of the sensing sector. After computing the coverage area for all feasible cases, each sensor selects the unit direction vector pointing toward the Voronoi vertex that yields the maximum coverage as its optimal orientation. If multiple sensors select the same vertex, the sensor achieving lower coverage reorients toward its next-best vertex to reduce overlap.

\subsection{Proposed orientation optimization algorithm}
We propose an orientation optimization algorithm that couples coverage maximization with robustness evaluation using the RRF, see Algorithm~\ref{algo:integrated}. The algorithm first computes the RRF for each sensor in the current configuration. If the minimum RRF across all sensors is at least $\rho_{\min}$, where $\rho_{\min}$ represents the minimum expected location error and $\rho_{\max}$ denotes a practical upper bound on possible deviations, the robust counterpart model \eqref{RC_model} is solved to optimize orientations under worst-case perturbations.

If the minimum RRF falls below $\rho_{\min}$, the algorithm switches to the robustified model \eqref{rrf-rc_model}, which explicitly enforces $\text{RRF} \geq \rho_{\min}$ for all sensors, ensuring resilience against admissible location deviations. This adaptive strategy applies additional robustness only when necessary, reducing conservatism and computational effort while maintaining reliable coverage. For simplicity, identical sensors with uniform sensing parameters are considered.
\begin{algorithm}[htbp]
\caption{Integrated Sensor Orientation Optimization with  RRF-based Robustness Evaluation}
\label{algo:integrated}
\begin{algorithmic}[1]
\Require Number of sensors $m$, initial positions $S = \{s_1, \dots, s_m\}$, sensing range bounds $r_s$ and $R_s$, maximum overlap $\epsilon$, minimum change in overall area $\delta$, minimum and maximum RRF threshold $\rho_{\min}$ and $\rho_{\max}$, sensing angle $\theta_H$
\Ensure Final optimal robust sensor orientations and corresponding maximum coverage

\State \textbf{Step 1: Initialization}
\State Set initial sensor positions $S$ and orientations
\State Compute initial Voronoi diagram for $S$ based on initial (nominal) sensor locations

\State \textbf{Step 2: Robustness Assessment}
\For{each sensor at $s_i,\; i=1:m$}
    \State Compute RRF $\rho_{i}$.
\EndFor

\State \textbf{Step 3: Voronoi-based Local Refinement}
\For{each sensor at $s_i \in S$}
    \State Identify Voronoi vertices for $s_i$ within allowable region
    \State For each vertex, compute coverage area at worst-case locations
    \State Choose the vertex maximizing average coverage of $s_i$ over its uncertainty set
\EndFor

\State \textbf{Step 4: Robust Coverage Optimization}
\For{each sensor at $s_i \in S$}
    \If{$\rho_{min} \geq \rho_{i} \text{ or } \rho_{i} \geq \rho_{\max}$}
        \State Sandwich $\rho_{i}$ it between RRF bounds.
    \EndIf
    \State Formulate and solve model \eqref{rrf-rc_model}.
    \State Record candidate robust orientations and coverage.
\EndFor

\State \textbf{Step 5: Cooperative Recalibration}
\Repeat
    \For{each pair $(s_i, s_j)$ sharing a Voronoi vertex}
        \If{overlap is more than threshold}
            \State Compare coverage areas $A_i$ and $A_j$
            \State Sensor with smaller coverage reorients to next-best vertex
        \EndIf
    \EndFor
\Until{No further beneficial adjustments possible (or overall change in area is lesser than $\delta$ threshold)}

\State \textbf{Output:} Final sensor orientations and total coverage.
\end{algorithmic}
\end{algorithm}

\section{Experimental analysis}  \label{sec:5}  
This section evaluates the proposed RRF-based framework by comparing nominal and perturbed deployments, analyzing Voronoi-cell-wise coverage, and studying the effect of sensing parameters on robustness and coverage. 

\subsection{Scenario setup} 
We consider a random deployment of $m$ aerial directional sensors in a $1000 \times 1000$ square region. Each sensor is associated with its Voronoi cell and contributes coverage based on its angular view $\theta_H$, with the objective of maximizing total coverage while accounting for uncertainty and overlap. Although we present results under a common value of $\theta_H$, the formulation is flexible and can be readily extended to accommodate heterogeneous angular ranges across different sensors, thereby enhancing the applicability of the analysis to more practical scenarios. All simulations and numerical experiments were implemented in Python.

Orientations yielding coverage below a threshold $\lambda$ are discarded; sensors with no valid orientation are randomized or put to sleep, thereby conserving energy and extending the overall network lifetime. An iterative process resolves excessive overlap using a threshold $\epsilon$ and terminates when coverage improvement falls below $\delta$.

\subsection{Case study: Voronoi-based coverage}
\begin{figure}
\centering
\begin{subfigure}{.48\textwidth}
  \centering
  \includegraphics[width=\linewidth]{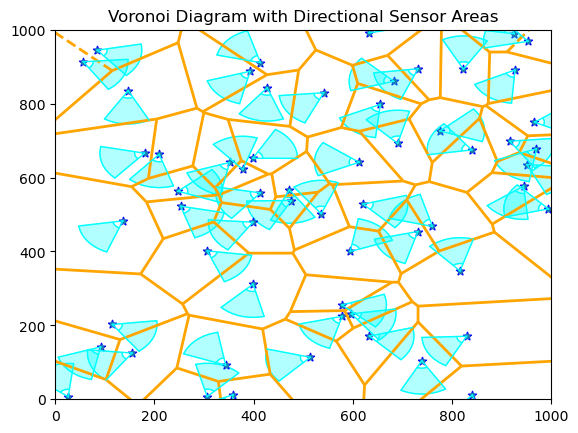}
  \caption{Best-case initial orientation}
  \label{Fig_4a}
\end{subfigure}
\hfill
\begin{subfigure}{.48\textwidth}
  \centering
  \includegraphics[width=\linewidth]{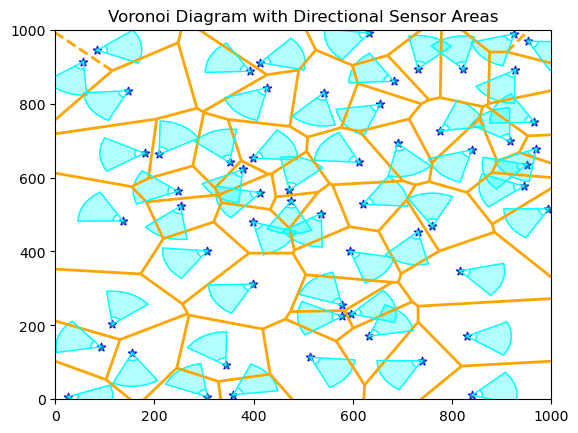}
  \caption{Best-case robustified orientation}
  \label{Fig_4b}
\end{subfigure}
\caption{Voronoi diagrams showing coverage under best-case (nominal) sensor deployment.}
\label{fig:best_case}
\end{figure}

To evaluate the proposed robustified orientation framework, we compare its performance under nominal (best-case) and perturbed (worst-case) sensor placements. Figures~\ref{Fig_4a} and \ref{Fig_4b} show that robustified orientations steer sensors toward Voronoi vertices that preserve coverage under perturbations, whereas nominal orientations often result in fragmented coverage, uncovered regions, and increased overlaps, reducing efficiency.

For idealized deployments with infinitely many sensors, convergence is reached when the incremental gain in total coverage becomes negligible. In practical scenarios with a finite number of sensors $m$, the algorithm terminates once all sensors are either assigned feasible orientations or placed in a non-active state, ensuring computational efficiency under realistic resource constraints.

\begin{figure}
\centering
\begin{subfigure}{.48\textwidth}
  \centering
  \includegraphics[width=\linewidth]{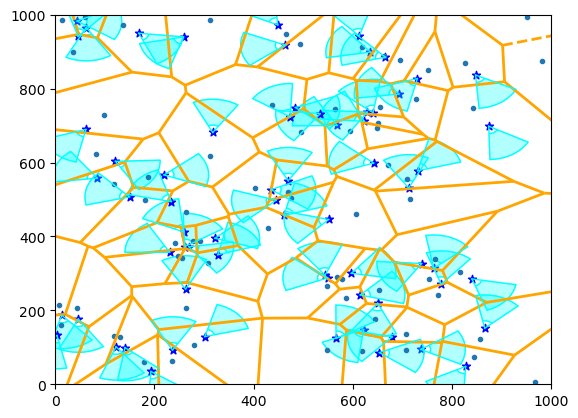}
  \caption{Worst-case initial orientation}
  \label{Fig_5a}
\end{subfigure}
\hfill
\begin{subfigure}{.48\textwidth}
  \centering
  \includegraphics[width=\linewidth]{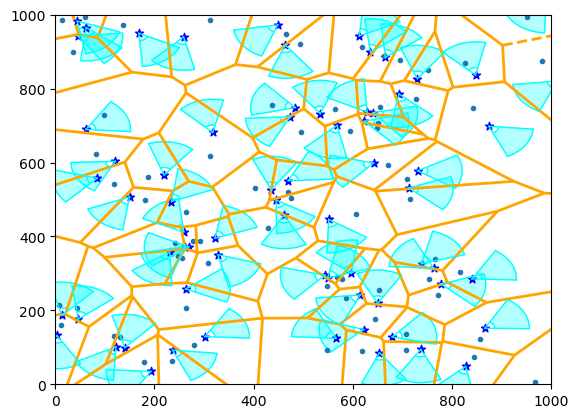}
  \caption{Worst-case robustified orientation}
  \label{Fig_5b}
\end{subfigure}
\caption{Voronoi diagrams showing coverage under worst-case (perturbed) sensor deployment.}
\label{fig:worst_case}
\end{figure}

To further assess performance under uncertainty, we perturb the same nominal sensor positions to model worst-case deployments. Figures~\ref{Fig_5a}–\ref{Fig_5b} compare (i) initial unadjusted orientations and (ii) robustified orientations, showing that the robustified model preserves coverage under significant perturbations while reducing fragmentation and overlap, whereas random orientations yield non-uniform and lower coverage. This advantage is more pronounced in the perturbed case, where robustified orientations clearly outperform the initial configuration. Over 1000 randomized deployments, the average effective coverage increases from $31{,}584.37$ (initial) to $120{,}437.89$ (robustified), demonstrating the robustness and scalability of the proposed framework in uncertain environments.

\subsection{Comparative coverage analysis}
\begin{figure}
    \centering
    \begin{subfigure}{0.48\textwidth}
        \centering
        \includegraphics[width=\linewidth]{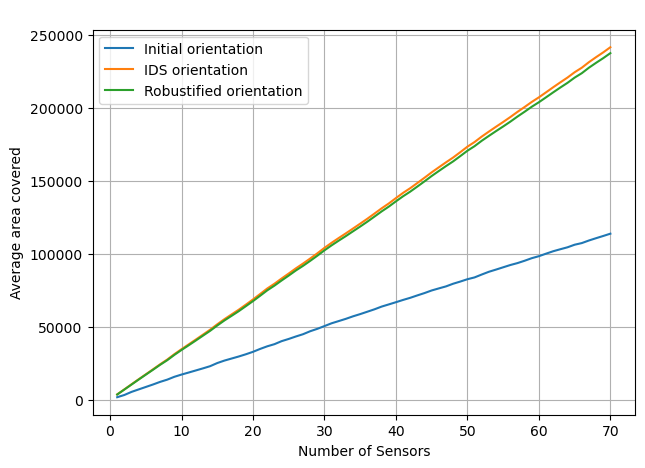}
        \caption{Coverage comparison under nominal deployment.}
        \label{fig:nominal_comparison}
    \end{subfigure}
    \hfill
    \begin{subfigure}{0.48\textwidth}
        \centering
        \includegraphics[width=\linewidth]{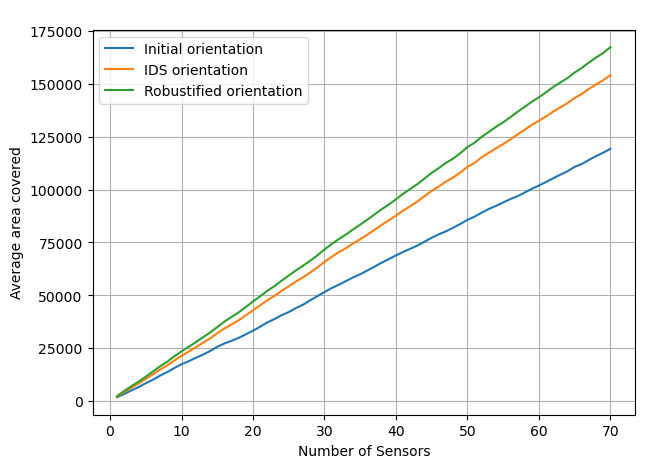}
        \caption{Coverage comparison under perturbed deployment.}
        \label{fig:worst_comparison}
    \end{subfigure}
    \caption{Coverage comparison (IDS, Random, and Robustified) for nominal (left) and perturbed (right) deployments.}
    \label{fig:coverage_comparison}
\end{figure}

We compare the proposed robustified orientation with two baselines: random orientation and IDS orientation \cite{sung2014voronoi}, under both nominal and perturbed deployments (Figs.~\ref{fig:coverage_comparison}). For exact sensor locations, robustified and IDS orientations achieve nearly identical coverage and both clearly outperform random orientation. The slightly lower coverage of the robustified method arises from its conservative optimization over the full uncertainty set rather than only nominal positions, reflecting a deliberate trade-off to ensure resilience against deployment perturbations.

\begin{figure}
    \centering
    \begin{subfigure}{0.48\textwidth}
        \centering
        \includegraphics[width=\linewidth]{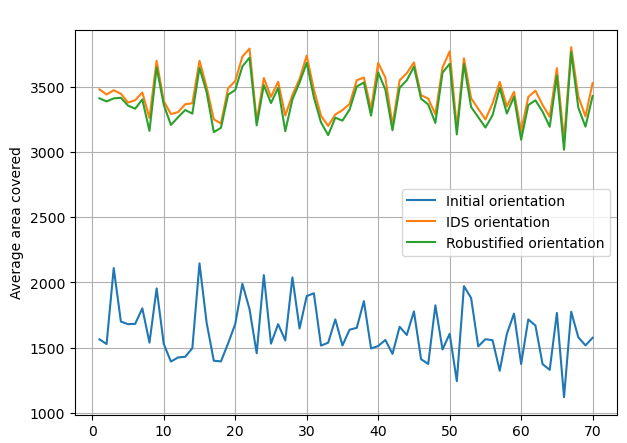}
        \caption{Cellwise coverage comparison under nominal deployment.}
        \label{fig:nominal_comparison}
    \end{subfigure}
    \hfill
    \begin{subfigure}{0.48\textwidth}
        \centering
        \includegraphics[width=\linewidth]{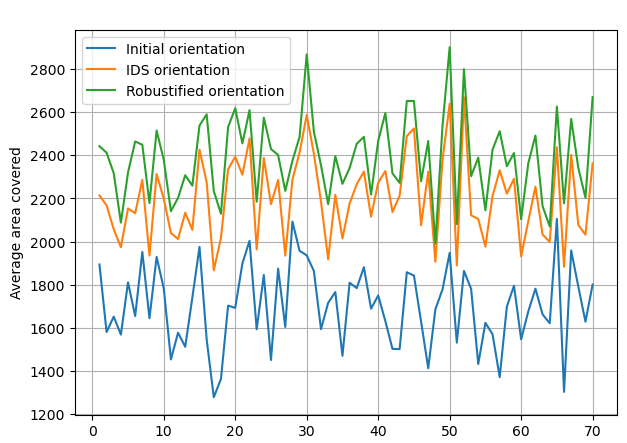}
        \caption{Cellwise coverage comparison under perturbed deployment.}
        \label{fig:worst_comparison}
    \end{subfigure}
    \caption{Voronoi-cell wise area distribution across sensors for different orientation strategies}
    \label{Fig_7}
\end{figure}

Under perturbed (worst-case) sensor placements, the performance gap between IDS and the robustified orientation becomes pronounced (Fig.~\ref{fig:worst_comparison}). IDS coverage degrades significantly under uncertainty, whereas the robustified approach consistently maintains larger effective coverage areas. Random orientation produces highly irregular and sparse coverage and performs so poorly that it is not a meaningful baseline. These results show that, in realistic deployments with unavoidable placement uncertainty, robustified orientation offers reliable coverage guarantees, even at the cost of a small reduction in nominal performance, which is often preferable to fragile best-case solutions. Figure~\ref{Fig_7} shows the Voronoi cell–wise coverage distribution for the three orientation strategies. 

\subsection{Parametric sensitivity analysis}
We study the sensitivity of the proposed framework by varying six key parameters one at a time, while keeping the others fixed, with the resulting robustified coverage summarized in Table~\ref{tab:parametric}. The nominal data used for this analysis is: 70 sensors with $\rho_{min}=25, \rho_{max}=35, \theta_H=60, r_s=30$ and $R_s=100$. This analysis evaluates how coverage responds to changes in sensing, deployment, and uncertainty settings, and confirms that the observed robustness trends remain consistent across different parameter configurations.

Wider horizontal sensing angles $\theta_H$ yield the largest coverage gains, though practical limits arise from hardware and energy constraints. Increasing the number of sensors improves coverage only up to a point, after which overlaps cause diminishing returns. This observation provides a valuable guideline for identifying an approximate upper bound on the number of sensors that should be deployed for cost-effective yet robust coverage. Regarding uncertainty, increasing $\rho_{\max}$ has little effect since the computed RRF remains smaller, whereas larger $\rho_{\min}$ values reduce coverage once robustness requirements become too restrictive. Finally, enlarging the maximum sensing radius increases coverage, while a larger minimum radius reduces it. These effects are consistently observed under both nominal and perturbed deployments, as reflected in Figure~\ref{Fig_7}, that depict Voronoi cell-wise area distribution across sensors for different orientation strategies.

\begin{table}[h!]
\centering
\scriptsize
\caption{Coverage comparison under varying parameters}
\label{tab:parametric}
\centering
\begin{tabular}{|p{2.2cm}|p{1cm}|p{1.5cm}|p{1.5cm}|p{1.5cm}|p{1.5cm}|}
\hline
\textbf{Parameter} & \textbf{Value Tested} & \textbf{Initial IDS}  & \textbf{Initial Robust}  & \textbf{Perturbed IDS}  & \textbf{Perturbed Robust} \\ \hline
\multirow{5}{*}{\shortstack{Number of\\ sensors $(m)$}} 
  & 50  & 197545.61  & 195675.73  & 147914.41  & 156159.29 \\ \cline{2-6}
  & 60  & 220417.59 & 217498.01 & 152840.40 & 163758.72 \\ \cline{2-6}
  & 70  & 241292.74 & 237235.40 & 154036.13 & 167248.34 \\ \cline{2-6}
  & 80  & 255771.48 & 251112.61 & 152637.65 & 165882.80 \\ \cline{2-6}
  & 90 & 268353.03 & 262454.90 & 149877.85 & 167063.18 \\ \hline
 
\multirow{6}{*}{\shortstack{Horizontal\\ angle $(\theta_H)$}} 
  & 30$^\circ$  &  133554.01 & 130680.25 & 89610.48 & 95551.78\\ \cline{2-6}
  & 60$^\circ$  & 241292.74 & 237235.40 & 154036.13 & 167248.34 \\ \cline{2-6}
  & 90$^\circ$  & 323679.37 & 318585.42 & 204926.30 & 222824.72 \\ \cline{2-6}
  & 120$^\circ$ & 400477.87 & 394865.70 & 257141.78 & 276700.57\\ \cline{2-6}
  & 150$^\circ$ & 456541.96 & 450345.93 & 293991.55 & 317964.21\\ \cline{2-6}
  & 180$^\circ$ & 509893.83 & 503138.02 & 334518.25 & 364027.01 \\ \cline{2-6}
  & 360$^\circ$ &  687955.95 & 687953.39 & 667221.44 & 667221.92 \\ \hline
  
\multirow{7}{*}{\shortstack{RRF bounds\\ $(\rho_{min} $ and $ \rho_{max})$}} 
  & 5-15 & 241936.41 & 241354.46  & 205024.89  & 208239.12 \\ \cline{2-6}
  & 15-25 & 241627.33 & 239754.49 & 179744.56 & 187668.15 \\ \cline{2-6}
  & 25-35 & 241292.74 & 237235.40 & 154036.13 & 167248.34 \\ \cline{2-6}
  & 35-45 & 242777.91 & 234572.34 & 131342.83 & 149069.48 \\ \cline{2-6}
  & 45-55 & 243183.01 & 226536.65 & 106442.46 & 130358.23 \\ \hline
  
\multirow{7}{*}{\shortstack{Sensing range\\ $(r_s \text{ and } R_{s})$}} 
  & 15-70   & 157394.44 & 155231.87 & 121440.72 & 132197.89\\ \cline{2-6}
  & 25-90   & 218921.17 & 215528.41 & 148463.48 & 160372.28\\ \cline{2-6}
  & 30-100   & 241292.74 & 237235.40 & 154036.13 & 167248.34 \\ \cline{2-6}
  & 35-110   & 258450.53 & 254055.79 & 156500.41 & 170078.14\\ \cline{2-6}
  & 40-120   & 263541.57 & 259491.68 & 150681.15 & 163818.09 \\ \cline{2-6}
  & 45-140   & 275266.13  & 270891.43 & 144920.78 & 158521.37\\ \hline
\end{tabular}
\end{table}

\subsection{Key insights from experiments}
The experiments show that the RRF-based robustified orientation consistently improves coverage, remains resilient under perturbations, and adapts well to sensing parameter variations, confirming the practical value of robustness in aerial directional sensor networks. In nominal deployments, initial orientation achieved $109{,}209.47$ square units, IDS increased coverage to $206{,}829.61$ square units, and the robustified approach achieved $201{,}815.78$ square units, slightly lower than IDS but with added robustness. Under worst-case perturbations, coverage dropped to $31{,}584.37$ square units for initial orientation, improved to $111{,}505.36$ square units with IDS, and further increased to $120{,}437.90$ square units with the robustified strategy. Thus, although IDS can marginally outperform robustified orientation in nominal settings, the robustified framework ensures higher and more reliable coverage under uncertainty, making it a more effective deployment strategy. While the current study focuses on a static and homogeneous setting, the proposed framework is not limited to these assumptions. Note that the framework can be extended to heterogeneous sensors by incorporating sensor-specific parameters (e.g., varying sensing ranges and FoVs) into the RRF formulation, while dynamic settings can be handled by updating the RRF over time using a sequential or time-varying optimization approach.





\section{Conclusion and future work}  \label{sec:6}

This article addresses the practical challenge of coverage degradation in aerial directional sensor networks caused by unavoidable location uncertainty. By modeling sensing through Voronoi-partitioned ground regions and explicitly incorporating the radius of robust feasibility, we develop an orientation strategy that guarantees coverage feasibility under bounded perturbations. The proposed robustified orientation algorithm operates using local geometric information, preserves feasibility by construction, and limits redundant overlaps. Extensive simulations confirm that the method achieves stable and reliable coverage under worst-case sensor displacements while remaining competitive in nominal conditions. These results demonstrate that embedding RRF directly into orientation design is an effective and practical approach for uncertainty-aware aerial-sensing deployments.\\
Future work includes real-time distributed implementations, integration with mobile sensing and handling compound uncertainties. Extending the framework to multi-objective settings (coverage, energy, latency) further enhances its applicability to broader uncertainty-aware sensing problems.



\begin{thebibliography}{99}

\bibitem{adachi2022cooperative} Adachi, T., Hayashi, N. Takai, S.: Cooperative target tracking by multiagent camera sensor networks via Gaussian process. IEEE Access, IEEE. \textbf{10}, 71717-71727 (2022)

\bibitem{Atmaji} Atmaja, B. T., Ihsannur, H., Arifianto, D.: Lab-scale vibration analysis dataset and baseline methods for machinery fault diagnosis with machine learning. J. Vib. Eng. Technol. \textbf{12} (2), 1991–2001 (2024)

\bibitem{ben2009robust}  Ben-Tal, A., El Ghaoui, L., Nemirovski, A.: Robust Optimization. Princeton University Press, Princeton (2009)

\bibitem{ben2000robust} Ben-Tal, A., Nemirovski, A.: Robust solutions of linear programming problems contaminated with uncertain data. Math. Program. \textbf{88}, 411-424 (2000)

\bibitem{ben2002robust} Ben-Tal, A., Nemirovski, A., Roos, C.: Robust solutions of uncertain quadratic and conic-quadratic problems. SIAM J. Optim. \textbf{13} (2), 535-560 (2002)

\bibitem{chen2020radius} Chen, J., Li, J., Li, X., Lv, Y., Yao, J.C.: Radius of robust feasibility of system of convex inequalities with uncertain data. J. Optim. Theory Appl. \textbf{184} (2), 384-399 (2020)

\bibitem{chuong2016robust} Chuong, T., Jeyakumar, V.: Robust global error bounds for uncertain linear inequality systems with applications. Linear Algebra and its Applications, Elsevier. \textbf{493}, 183-205 (2016)

\bibitem{chuong2017exact} Chuong, T.D., Jeyakumar, V.: An exact formula for radius of robust feasibility of uncertain linear programs. J. Optim. Theory Appl. \textbf{173} (1), 203-226 (2017)

\bibitem{2} Dai, X., Hu, Y., Cui, D., Chai, T.: A disturbance decoupling generalized proportional-integral observer design for robust sensor fault detection. IEEE Transactions on Industrial Electronics. \textbf{70} (6), (2022)

\bibitem{Dasig} Dasig, D. D.: Implementing IoT and wireless sensor networks for precision agriculture. In Internet of Things and Analytics for Agriculture. Singapore: Springer Singapore. \textbf{2}, 23-44 (2019)

\bibitem{vanshika} Datta, V., Nahak, C.: A Radius of Robust Feasibility Approach to Directional Sensors in Uncertain Terrain. arXiv preprint arXiv:2510.19407 (2025)

\bibitem{eren2025simulation} Eren, O., Altin K.A.: A simulation-informed robust optimization framework for the design of energy efficient underwater sensor networks. Ad Hoc Networks, Elsevier. \textbf{178}, 103933 (2025)

\bibitem{gabrel2014recent} Gabrel, V., Murat, C., Thiele, A.: Recent advances in robust optimization: An overview. European J. Oper. Res. \textbf{235} (3), 471-483 (2014)

\bibitem{goberna2014robust} Goberna, M.A., Jeyakumar, V., Li, G., Vicente-P{\'e}rez, J.: Robust solutions of multiobjective linear semi-infinite programs under constraint data uncertainty. SIAM J. Optim. \textbf{24} (3), 1402-1419 (2014)

\bibitem{goberna2022radius} Goberna, M.A., Jeyakumar, V., Li, G., Vicente-P{\'e}rez, J.: The radius of robust feasibility of uncertain mathematical programs: A Survey and Recent Developments. European J. Oper. Res. \textbf{296} (3), 749-763 (2022)

\bibitem{goberna2015robust} Goberna, M.A., Jeyakumar, V., Li, G., Vicente-P{\'e}rez, J.: Robust solutions to multi-objective linear programs with uncertain data. European J. Oper. Res. \textbf{242} (3), 730-743 (2015)

\bibitem{goberna2016radius} Goberna, M.A., Jeyakumar, V., Li, G., Linh, N.: Radius of robust feasibility formulas for classes of convex programs with uncertain polynomial constraints. Oper. Res. Lett. \textbf{44} (1), 67-73 (2016)

\bibitem{goberna2021calculating} Goberna, M. A., Jeyakumar, V., Li, G.: Calculating radius of robust feasibility of uncertain linear conic programs via semi-definite programs. J. Optim. Theory Appl. Springer. \textbf{189} (2), 597-622 (2021)

\bibitem{goberna2018guaranteeing} Goberna, M. A., Jeyakumar, V., Li, G., P{\'e}rez, J. V.: Guaranteeing highly robust weakly efficient solutions for uncertain multi-objective convex programs. European Journal of Operational Research, Elsevier. \textbf{270} (1), 40-50 (2018)

\bibitem{[12]} Javaid, S., Fahim, H., Zeadally, S., He, B.: Self-powered sensors: Applications, challenges, and solutions. IEEE Sensors J. \textbf{23} (18), 20483–20509 (2023)

\bibitem{Khan} Khan, R., Saeed, U., Koo, I.: FedLSTM: A federated learning framework for sensor fault detection in wireless sensor networks. Electronics. \textbf{13} (24), 4907 (2024)

\bibitem{3} Khan, R., Saeed, U., Koo, I.: Robust Sensor Fault Detection in Wireless Sensor Networks Using a Hybrid Conditional Generative Adversarial Networks and Convolutional Autoencoder. IEEE Sensors Journal. \textbf{25} (8), 1558-1748 (2025)

\bibitem{Krishnamoorthy} Krishnamoorthy, S., Dua, A., Gupta, S.: Role of emerging technologies in future IoT-driven healthcare 4.0 technologies: A survey, current challenges and future directions. J. Ambient Intell. Humanized Comput. \textbf{14} (1), 361–407 (2023)

\bibitem{kouvelis2013robust} Kouvelis, P., Yu, G.: Robust Discrete Optimization and Its Applications. Springer, Berlin. \textbf{14} (2013)

\bibitem{li2021coverage} Li, W., Wang, X., Han, S.: Coverage enhance in boundary deployed camera sensor networks for airport surface surveillance. IEEE Access. \textbf{9}, 145728-145738 (2021) 

\bibitem{li2022uav} Li, Chen, H.: UAV enhanced target-barrier coverage algorithm for wireless sensor networks based on reinforcement learning. Sensors, MDPI. \textbf{22} (17), 6381 (2022)

\bibitem{liang2025enhanced} Liang, J. M., Mishra, S., Lin, C. H.: Enhanced PTZ Camera Dispatch Scheme for 3D Environments Based on Deep Reinforcement Learning. IEEE Transactions on Instrumentation and Measurement, IEEE. \textbf{74}, 1557-9662 (2025)

\bibitem{liers2022radius} Liers, F., Schewe, L., Th{\"u}rauf, J.: Radius of robust feasibility for mixed-integer problems. INFORMS Journal on Computing, INFORMS. \textbf{34} (1), 243-261 (2022)

\bibitem{luo2020maximizing} Luo, C., Hong, Y., Li, D., Wang, Y., Chen, W., Hu, Q.:Maximizing network lifetime using coverage sets scheduling in wireless sensor networks. Ad Hoc Networks, Elsevier. \textbf{98}, 102037 (2020)

\bibitem{mahfouz2023novel} Mahfouz, A. M., Ismail, A. S., El Sobky, W. I., Nasry, H.: A novel model for representing a plane target and finding the worst-case coverage in wireless sensor network based on Clifford algebra. EURASIP Journal on Wireless Communications and Networking, Springer. \textbf{2023} (1), 95 (2023)

\bibitem{Pawar} Pawar, P., Kumar, M. T.,Vittal K. P.: An IoT based intelligent smart energy management system with accurate forecasting and load strategy for renewable generation. Measurement. \textbf{152}, 107187 (2020)

\bibitem{ridolfi2023radius} Ridolfi, A.B., Vera de Serio, V.N.: A Radius of Robust Feasibility for Uncertain Farthest Voronoi Cells. Set-Valued Var. Anal. \textbf{31} (1) (2023)

\bibitem{Saeed2022} Saeed, U., Shah, S. Y., Ahmad, J., Imran, M. A., Abbasi, Q. H., Shah, S. A.: Machine learning empowered COVID-19 patient monitoring using non-contact sensing: An extensive review. J. Pharmaceutical Anal. \textbf{12} (2), 193–204 (2022)

\bibitem{sangwan2015survey} Sangwan, A., Singh, R. P.:Survey on coverage problems in wireless sensor networks. Wireless Personal Communications, Springer. \textbf{80}, 1475-1500 (2015)

\bibitem{Sobin} Sobin, C. C.: A survey on architecture, protocols and challenges in IoT. Wireless Pers. Commun.\textbf{112} (3), 1383–1429 (2020)

\bibitem{soyster1973convex} Soyster, A.L.: Convex programming with set-inclusive constraints and applications to inexact linear programming. Oper. Res. \textbf{21} (5), 1154-1157 (1973)

\bibitem{sung2014voronoi} Sung, T. W., Yang, C. S.: Voronoi-based coverage improvement approach for wireless directional sensor networks. J. Netw. Comput. Appl. \textbf{39}, 202-213 (2014) 

\bibitem{wang2017coverage} Wang, Y., Wu, S., Chen, Z., Gao, X., Chen, G.: Coverage problem with uncertain properties in wireless sensor networks: A survey. Computer Networks, Elsevier \textbf{123}, 200--232 (2017)

\bibitem{[13]} Xing, L.: Cascading failures in Internet of Things: Review and perspectives on reliability and resilience. IEEE Internet Things J. \textbf{8} (1), 44–64 (2021)

\bibitem{ye2008robust} Ye, W., Ordonez, F.: Robust optimization models for energy-limited wireless sensor networks under distance uncertainty. IEEE transactions on wireless communications, IEEE. \textbf{7} (6), 2161-2169 (2008)

\bibitem{1} Yuen, K. V., Hao, X. H., Kuok, S. C.: Robust sensor placement for structural identification. Structural Control and Health Monitoring. \textbf{29} (1), (2022)

\end{thebibliography}
\end{document}